 

\baselineskip=14pt
\parskip=10pt
\def\Tilde{\char126\relax}
\def\halmos{\hbox{\vrule height0.15cm width0.01cm\vbox{\hrule height
 0.01cm width0.2cm \vskip0.15cm \hrule height 0.01cm width0.2cm}\vrule
 height0.15cm width 0.01cm}}
\font\eightrm=cmr8  
\font\eighttt=cmtt8
\magnification=\magstephalf

\parindent=0pt
\overfullrule=0in
 
\bf
\centerline
{PROOF of CONWAY'S LOST COSMOLOGICAL THEOREM}
\rm
\bigskip
\centerline{ {\it Shalosh B. Ekhad{$^1$}
and Doron Zeilberger}\footnote{$^1$}
{\eightrm  \raggedright
\baselineskip=9pt
Department of Mathematics, Temple University,
Philadelphia, PA 19122, USA. 
{\eighttt [ekhad,zeilberg]@math.temple.edu}
{\eighttt http://www.math.temple.edu/\Tilde [ekhad,zeilberg]   .}
Supported in part by the NSF.
First version: May 1, 1997. 
This version: July 9, 1997. 
} 
}
 
One of the most intriguing sequences  
([CG][F][SP][V])
is Conway's[C]
{\bf 1, 11, 21, 1211, 111221, 312211, 13112221, 1113213211,  ...}.
It is defined by the rule $C_0:=1$, and
$C_i:=JHC(C_{i-1})$, for $i>0$, where $JHC$ is Conway's
audioactive operator:
$$
JHC(a_1^{m_1} a_2^{m_2} \dots a_r^{m_r}):=
m_1 a_1 m_2 a_2  \dots  m_r a_r  \quad .
$$
Here $a^m$ is shorthand for ``$a$ repeated $m$ times''
(and we agree that the description is
optimal, i.e. $a_{i} \neq a_{i+1}$.)
We assume familiarity with Conway's charming article[C].
 
Conway proved that his sequence has the property
$length(C_{i+1})/length(C_i) \rightarrow \lambda$, where
$\lambda=1.303577269..$ is {\it Conway's constant}.
He also stated that, more generally, if one starts
with an {\it arbitrary} non-empty finite string of integers, $B_0$, 
(except `boring old 22'), and
defines $B_{i}:=JHC(B_{i-1})$, $i>0$, then still
$length(B_{i+1})/length(B_i) \rightarrow \lambda$.
This is an immediate consequence of
 
{\bf The Cosmological Theorem:} 
There exists an integer $N$ such
that every string decays in at most $N$ days to 
a compound of common and transuranic elements.
 
Conway stated
that two independent proofs {\it used to exist}, one by himself
and Richard Parker (that only proved that $N$ existed), 
and another one by Mike Guy (that actually proved that  
one may take $N=24$,
and that it was best possible). Unfortunately
both proofs were lost. Here we announce a new proof
(which establishes that one may take $N=29$; with more computations
one should be able to rederive (or else refute) Guy's sharp
$N=24$).
 
The Cosmological Theorem is an immediate consequence of the following lemma.
 
{\bf Lemma}: The length of any atom in 
the splitting of a $9$-day-old string is $\leq 80$. Every
such atom decays, in at most $20$ days, into stable or transuranic elements.
 
The lemma is {\it proved} by
typing {\tt Cosmo(8);} in the Maple package {\tt HORTON},
accompanying this announcement. The procedure {\tt Cosmo}
computes iteratively all non-splittable
strings of length $i$ ($i=1,2, \dots$) that might conceivably
be substrings (`chunks') of an atom in the
splitting of a $9$-day-old-string
(by backtracking, examining its possible ancestors
up to (at most) $8$ days back and rejecting those that 
lead to grammatically incorrect ancestors,
see examples below). 
Every time a string of length $i$
is accepted, its longevity (number of days it takes to decay
to stable or transuranic elements) is computed, and checked whether
it is finite. The maximal longevity turned out to be $20$.
The program halts if and when an $i$ is reached for which the set
of such conceivable strings of length $i$ is empty.
 
If the program halts (it did for us),
then the Lemma, and hence the Cosmological Theorem, are proved.
In fact it halted after $i=80$, implying that there do not
exist atoms of length $>80$ that occur in 
the splitting of {\it mature}
(i.e. $9$-day-old) strings, and that all the atoms have bounded
($\leq 20$, as it turned out) longevity. 
We also get that the longevity of an arbitrary
string is $\leq 9+20=29$.
The input and output files may be obtained from our websites.\halmos
 
The Maple package {\tt HORTON}, available from the authors'
websites, also rederives many other results
in Conway's paper, in particular it finds all the stable elements
{\it ab initio}, finds the minimal
polynomial for $\lambda$, finds the abundance of
all the stable elements,
and computes the longevity of any string. 
We refer the reader to the on-line documentation and to the source code.
 
{\bf Details}
 
Recall  that Conway proved that it suffices to consider strings
on $\{ 1,2,3\}$.
Let's call a chunk that starts with a comma {\it female}, and a chunk
that does not, {\it male}. Any chunk could have come either
from a father or a mother
(but of course not from both). Define
$$
ParentOfGirl( a_1 a_2 a_3 a_4 a_5 a_6  \dots):=
a_2^{a_1} a_4^{a_3} a_6^{a_5} \dots \quad ,
$$
and
$$
ParentOfBoy( a_1 a_2 a_3 a_4 a_5 a_6  \dots):=
a_1 a_3^{a_2} a_5^{a_4} \dots \quad .
$$
Since a parent of a chunk may be either female or male
(but we have no way of knowing), any chunk has two potential
parents (but of course only one actual one), (up to) four 
(potential) grandparents
(some of them may disqualify on the grounds of being grammatically
incorrect), and so on. Now there are lots of chunks that can't
possibly be factors of a mature string. Take for
example the female chunk ``$,12, 32,$''. It can't be
a chunk of a $1$-day-old string (as observed in [C]), since
starting at day $1$, all strings are descriptive, and
``$,12, 32,$'' would have been abbreviated  ``$,42,$''. So
we can {\it eliminate from the outset} any female string
of the form ``$, a c, b c,$'', in strings that are older than $0$ days.
 
Now consider the female string ``$,32,33,$''. It 
is grammatically correct, and so can be a chunk of a $1-$day-old string.
Her parent is:
``$222333$''. If the parent is a father, then it
is punctuated ``$2, 22, 33,3$'', which is grammatically incorrect,
and if the parent is a mother, then it would be
``$,22,23,33,$'', that is equally grammatically incorrect.
Hence we can conclude that ``$,32,33,$'', while it may
be a chunk of an atom in the splitting
of a $1$-day-old string, can't possibly be
such a chunk of a $2$-day-old string.
 
Consider on the other hand the female string ``$,12,21,$''. Her
father is ``$2,11$'', and her mother is ``$,21,$''. Her paternal
grandparent is $21$ who is OK, being two-lettered.
Hence we can't rule out ``$,12,21,$'' as a possible chunk in
an $L$-day old string (for any $L>0$).
 
Let's define $U_L(i)$ as the set of female strings 
on the alphabet $\{1,2,3\}$ of length $2i$ that
do not split, and that 
have at least one
$(great)^{L-2}$ grandparent who is grammatically
correct, or some $(great)^{j-2}$ grandparent ($j<L$) that is two-lettered.
In order to find $U_L(i)$, we take all the survivors
that made it to $U_L(i-1)$, append all the nine possible
endings $11,12,13, 21,22,23,31,32,33$, and see which of the resulting 
female strings 
do not split and in addition
survive up-to-depth-$L$ genealogical screening.
Whenever we induct a new
member to $U_L(i)$, we also compute her longevity, 
and the longevity of her male extensions of length $2i+1$:
``$1,w$'',``$2,w$'',``$3,w$'',
her female extensions of length $2i+1$:
``$w,1$'',``$w,2$'',``$w,3$'', and her nine male extensions of
length $2i+2$: 
``$1,w,1$'', ``$1,w,2$'', ``$1,w,3$'', ``$2,w,1$'', ``$2,w,2$'', 
``$2,w,3$'', ``$3,w,1$'', ``$3,w,2$'', ``$3,w,3$''.
We always keep track of the maximum longevity to-date.
 
If for some $L$ ($L=8$ worked), and some $i$ 
(we got $i=40$), $U_L(i)=\emptyset$, 
and the longevity record, $M$, is finite (we got $M=20$), then
it follows that any string $w$ decays into stable 
or transuranic elements in at most $M+L+1$ days.
Indeed, $w':=JHC^{9}(w)$ is a $9$-day-old string.
Split $w'$ into atoms. We know that each of these atoms
has length
$\leq 80$. Because if such an atom were a female of even length $\geq 82$,
then her head, consisting of the first $80$ letters, would have
been a member of $U_8(40)$, contradicting the fact that this
set is empty. Similarly if such an atom were a female of odd length, or
a male, then an appropriate factor would be a member of the empty $U_8(40)$.
So every atom in the splitting of $w'$ has length $\leq 80$.
Furthermore, all these atoms either belong to $U_L(i)$ for some $i<40$, 
or have one of the forms
$1v$, $2v$, $3v$, $v1$, $v2$, $v3$, 
$1v1$, $1v2$, $1v3$, $2v1$, $2v2$, $2v3$, $3v1$, $3v2$, $3v3$, 
where
$v \in U_L(i)$, for some $i<40$. But all these strings
were tested for finite longevity by {\tt Cosmo}, and turned out
to have longevity $\leq 20$. Since 
each of the atoms in the splitting of $w'$ decays in at most
$20$ days, so does $w'$,
and hence $w$ decays in at most $20+8+1=29$ days. \halmos
 
{\bf On A Posteriori Trivial Theorems: The Ultimate
Proof of the Four-Color Theorem Should Emulate our Proof}
 
Some statements are {\it a priori} trivial (regardless of whether there are
true or false), for example that there do not exist
projective planes of order $10$ (proved by Clement Lam), or the still open
conjecture that White can always win at Chess.
Of course by {\it trivial}
we mean {\it modulo a finite amount of calculations}.
 
Other statements are only {\it a posteriori} trivial.
Many conjectures seem only to be a posteriori trivial if
they turn out to be false, and are rendered trivial by
exhibiting a counter-example.
For example Euler's 
conjecture that $A^4+B^4+ C^4=D^4$ is unsolvable,
disproved by Elkies. Of course, in the
G\"odelian sense, every decidable statement is
a posteriori trivial. Both proof and disproof,
being finite, could be eventually found by exhaustive search.
 
However, the Cosmological Theorem turned out to be
{\it a posteriori trivial} in a more genuine `object-oriented'
sense. We wrote a program that iteratively computes
$U_L(i)$, and a proof would be obtained if $U_L(i)$ 
is empty for some $L$ and $i$. A priori, we have
no way of knowing whether we would be successful.
If our civilization would die, or more realistically, the
program would run out of memory, we would never know whether
it was `never' or `not yet'. All we can do is hope.
Also, suppose that the program does not halt.
In that case it does not mean that
the statement of the conjecture is false. All it means is that
our particular approach failed. So you only win if and when the
proof-program halts. 
But, once that happens,
in order to check the {\it validity of our  proof}, it is a
waste of time to examine the members of $U_8(i)$ for $i=1,2,3, \dots$,
and to examine the decay process for each accepted chunk. All we
need is to {\it check the computer program}. Once the logic of the
program has been checked, all we have to do is, bet on an $L$,
say $L=8$, and type {\tt Cosmo(8);}, and hope
that it would halt in our lifetime. If it does, that's all there
is to it, and we have a one-line proof.
 
A celebrated example of an {\it a posteriori trivial} theorem
is Appel and Haken's Four-Color-Theorem. Their approach was
to find a finite unavoidable set of 
reducible configurations. 
The original proof[AHK] used
an excessive amount of human effort. This was considerably
improved in the new proof by Robertson, Sanders, Seymour, and Thomas[RSST],
but this is not the {\it ultimate proof}. Eventually
one should be able to type {\tt Prove4CT();}, and the
truth of the theorem should be implied by the halting of
the program. In order to check the validity, the checker would not
need to see any specific configuration. Everything should be
done internally and silently by the computer. All that the checker
would have to do is {\it check the program}. 
 
And who knows? Perhaps the non-existence of odd perfect numbers,
the $3n+1$ problem,
the Riemann Hypothesis, etc. etc. are all {\it a posteriori} trivial?
Let's hope that their proof-programs will halt in our lifetime.

{\bf Acknowledgement.} We wish to thank John Conway, Steve Finch, and
Ilan Vardi for helpful remarks on an earlier version.
 
{\bf References}
 
[AHK] K. Appel, W. Haken, and J. Koch, {\it Every planar
map is four-colorable},  Illinois J. Math. {\bf 21}(1977),
429-567.
 
[C] J.H. Conway, {\it The weird and wonderful chemistry of
audioactive decay}, in: ``{\it Open Problems in Communication
and Computation}'', T.M. Cover and B. Gopinath, eds., Springer,
1987, pp. 173-188.
 
[CG] J.H. Conway and R.K. Guy, {\it ``The Book of Numbers''},
Copernicus, Springer, 1996.
 
[F] S. Finch, {\it ``Favorite Mathematical Constants Website''},
{\hfill \break
\tt http://www.mathsoft.com/asolve/constant/cnwy/cnwy.html}.
 
[RSST] N. Robertson, D.P. Sanders, P. Seymour, and R. Thomas,
{\it A new proof of The Four-Color Theorem}, Elect. Res.
Announc. of the Amer. Math. Soc. {\bf 2}(1996), no. 1, 17-25 (electronic).
 
[SP] N.J.A. Sloane and S. Plouffe, {\it ``The Encyclopedia
of Integer Sequences''}, Academic Press, 1995.
 
[V] I. Vardi, {\it ``Computational Recreations in Mathematica''},
Addison-Wesley, 1991.

\bye